\documentclass{proc-l}
\pdfoutput=1
\issueinfo{130}{09}{September}{2002}
\pagespan{2525}{2534}
\dateposted{February 4, 2002}
\PII{S 0002-9939(02)06366-9}
\copyrightinfo{2002}{American Mathematical Society}
\usepackage{amscd}
\usepackage{euscript}
\usepackage{hyperref}

\renewcommand{\MR}[2]{MR \href{http://www.ams.org/mathscinet-getitem?mr=#1:#2}{{\bf #1}:{#2}}}

\newtheorem{theorem}{Theorem}[section]
\newtheorem{lemma}[theorem]{Lemma}

\theoremstyle{definition}

\newtheorem*{algorithmn}{Algorithm}
\newtheorem{algorithm}{Algorithm}
\newtheorem*{inversealgorithm}{Inverse Algorithm}

\theoremstyle{remark}
\newtheorem{remark}[theorem]{Remark}

\numberwithin{equation}{section}

\begin{document}
\title[A PROOF OF PIERI'S FORMULA]{A proof of Pieri's formula using the 
generalized Schensted insertion algorithm for rc-graphs}

\author{Mikhail Kogan}
\address{Department of Mathematics, Northeastern University, Boston,
Massachusetts 02115}
\email{misha@research.neu.edu}

\author{Abhinav Kumar}
\address{Department of Mathematics, Massachusetts Institute of Technology,
Cambridge,  Massachusetts 02139}
\email{abhinavk@mit.edu}

\commby{John R. Stembridge}
\date{November 17, 2000 and, in revised form, April 6, 2001}
\begin{abstract}
We provide a generalization of the Schensted insertion algorithm for
rc-graphs of Bergeron and Billey.  The new algorithm is used to
give a new proof of Pieri's formula.
\end{abstract}

\maketitle

\section{Introduction}

RC-graphs were first defined by Fomin and Kirillov \cite{fk} and later
studied by Bergeron and Billey \cite{bb}. They encode the monomials
contributing to the expansion of Schubert polynomials. These polynomials
were introduced by Lascoux and Schutzenberger \cite{ls1}, \cite{ls2} and
described at length by Macdonald \cite{m} and by Manivel~\cite{ma}.

A central problem in the theory of Schubert polynomials is to provide
effective ways of computing the generalized Littlewood-Richardson
coefficients $c^{u}_{vw}$ in the expansion
\begin{equation}
P_vP_w = \sum_u c^{u}_{vw} P_u,
\end{equation}
where $P_w$ is the Schubert polynomial of a permutation $w$.

The first attempt to compute the generalized Littlewood-Richardson
coefficients using rc-graphs was made in \cite {bb}, where Monk's formula
was proved using a generalized Schensted insertion algorithm. Later, in
\cite{k}, this generalized algorithm was used to compute a more general
set of Littlewood-Richardson coefficients. Unfortunately, this algorithm
does not work in general and, in particular, fails to give a proof
of Pieri's formula. In this paper we will modify this algorithm to
inserting a whole row of elements at once instead of an element by element
insertion. This generalization will prove Pieri's formula.

Let us give a statement of Pieri's formula, which was originally proved in
\cite{ls1} and later reproved by other methods \cite{s}, \cite{p},
\cite{w}. Our way of stating Pieri's formula is not standard, but in
Section 4 we show how to deduce the standard formulation of Pieri's formula
from Theorem \ref{pieri}.

\setcounter{section}{4}
\begin{theorem}
\[
P_w P_{\sigma[r,m]} = \sum P_{w'},
\]
where the sum is over all $w' = wt_{a_1,b_1}...t_{a_m,b_m}$ such that the
$b_i$'s are distinct and greater than $r$, $a_1 \leq a_2 \leq ... \leq
a_m\leq r$, satisfying $l(w')= l(w)+m$ and $w'(b_i) < w'(b_j)<
w'(a_i)$ for every $i<j$ with $a_i=a_j$.
\end{theorem}
\setcounter{section}{1}

In the above statement, $\sigma[r,m]$ stands for the permutation
$[1,2,...,r-1, r+m, r, r+1,...]$. Note that $P_{\sigma[r,m]}$ is equal
to the homogeneous symmetric polynomial $h_m(x_1,..., x_r)$.

We will define Schubert polynomials and rc-graphs, and point out their
relationship, in Section 2. The generalized insertion algorithm will be
given in Section 3, while in Section 4 we will prove Pieri's formula by
providing the inverse algorithm.

\section{RC-graphs}

Let us decide on some notation and state some basic facts about
permutations and Schubert polynomials first. For a permutation $w \in
S_n$, where $S_n$ is the symmetric group on $n$ elements, we can write $w$
as $[w(1),...,w(n)]$. We let $t_{ab}$ denote  the
transposition which takes $a$ to $b$ and vice versa, leaving other
elements fixed, and write $s_i$ for $t_{i,i+1}$. It can be shown that
$s_1,...,s_{n-1}$ generate $S_n$ with the relations
\begin{eqnarray*}
s_i^2 &=& 1, \cr
s_i s_j &=& s_j s_i \textrm{   if } |i-j| \geq 2, \cr
s_is_{i+1}s_i &=& s_{i+1}s_is_{i+1}.
\end{eqnarray*}

For a permutation $w \in S_n$, we let the length $l(w)$ of $w$ be the
number of inversions of $w$, i.e. the number of pairs $(i,j)$ with $i<j$
and $w(i) > w(j)$. A string $a_1a_2...a_p$ such that $s_{a_1}...s_{a_p} =
w$ with $p$ minimal is called a reduced word for $w$. It can be shown
that this $p$ is just the length of $w$, and any two reduced word
decompositions can be transformed into each other using the last two
relations given above. Denote by $R(w)$ the set of all the reduced words
for $w$.

This enables us to define Schubert polynomials. Let $\partial_i$ be an operator
which, acting on a polynomial $f(x_1,...,x_n)$ to the left, yields the
polynomial\begin{displaymath}
\frac{f(x_1,...,x_i,x_{i+1},....,x_n) - f(x_1,...,x_{i+1},x_i,...,
x_n)}{x_i-x_{i+1}}.
\end{displaymath}

It can be verified that
\begin{eqnarray*}
\partial_i^2 &=& 0, \cr
\partial_i \partial_j &=& \partial_j \partial_i \textrm{  if } |i-j| \geq 2 ,\cr
\partial_i \partial_{i+1} \partial_i &=& \partial_{i+1} \partial_i
\partial_{i+1}.
\end{eqnarray*}

Thus for any permutation $w$, we may define $\partial_w$ as $\partial_{a_1}...\partial_{a_p}$, where $a_1...a_p$ is {\em any} reduced word for $w$.
Now the Schubert polynomial $P_w$ is defined to be
\begin{eqnarray*}
P_w = \partial_{w^{-1}w_0} x_1^{n-1}x_2^{n-2}...x_{n-1}^1x_n^0,
\end{eqnarray*} where $w_0 = [n,...,1]$ is the longest permutation in
$S_n$.

Now we come to rc-graphs. Given a reduced word $a_1...a_p$ for $w$, we say
that a sequence $\alpha_1,...,\alpha_p$ is compatible for this word if
\begin{eqnarray*}
\alpha_1 \leq \alpha_2 \leq ...\leq \alpha_p  ,\cr
\alpha_i \leq a_i,  \forall \; i ,\cr
\alpha_i < \alpha_{i+1} \textrm{  if } a_i < a_{i+1}.
\end{eqnarray*}
Denote by $C({\bf a})$ the set of all compatible sequences for ${\bf a}\in
R(w)$.

Given such a compatible sequence, we make an rc-graph by placing
intersections (see Figure 1) at all the positions $(\alpha_k, a_k-\alpha_k
+1)$; these form a subset of $\{1,2,...\} \times \{1,2,...\}$. Now we draw
strands starting at each row and winding their way up. Wherever there is
an intersection symbol, we make the two strands there intersect. It is
easy to show by induction that the strand starting at the $i^{th}$ row
ends up at the $w(i)^{th}$ column and that no two strands intersect twice.
Conversely, a graph with intersections, such that no two strands intersect
twice, immediately gives us a reduced word for the permutation and also a
compatible sequence (we look at the intersections in order of increasing
row number; going from right to left in each row, the row number gives the
$\alpha_i$, and to get $a_i$ we add the row number and column number and
subtract 1). For more details see \cite{bb}.

\setlength{\textheight}{8.1in}

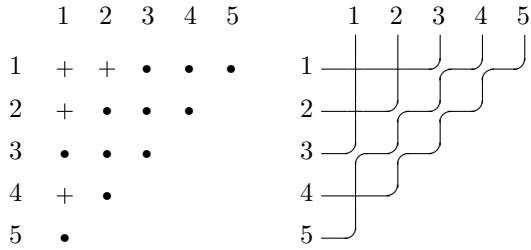
\begin{figure}
\vspace{-10pt}
\begin{picture}(200,134)
\put (35,68){\circle*{3}}
\put (35,36){\circle*{3}}
\put (51,84){\circle*{3}}
\put (51,68){\circle*{3}}
\put (51,52){\circle*{3}}
\put (66,100){\circle*{3}}
\put (66,84){\circle*{3}}
\put (66,68){\circle*{3}}
\put (82,100){\circle*{3}}
\put (82,84){\circle*{3}}
\put (98,100){\circle*{3}}

\put (32,84){\line(1,0){6}}
\put (35,81){\line(0,1){6}}

\put (32,100){\line(1,0){6}}
\put (35,97){\line(0,1){6}}

\put (32,52){\line(1,0){6}}
\put (35,49){\line(0,1){6}}

\put (48,100){\line(1,0){6}}
\put (51,97){\line(0,1){6}}

\put(14,98){1}
\put(14,82){2}
\put(14,66){3}
\put(14,50){4}
\put(14,34){5}

\put(32,117){1}
\put(48,117){2}
\put(64,117){3}
\put(80,117){4}
\put(96,117){5}

\put (132,100){\line(1,0){10}}
\put (132,84){\line(1,0){10}}
\put (132,68){\line(1,0){10}}
\put (132,52){\line(1,0){10}}
\put (132,36){\line(1,0){10}}
\put (148,100){\line(1,0){10}}
\put (148,84){\line(1,0){10}}
\put (148,68){\line(1,0){10}}
\put (148,52){\line(1,0){10}}
\put (164,100){\line(1,0){10}}
\put (164,84){\line(1,0){10}}
\put (164,68){\line(1,0){10}}
\put (180,100){\line(1,0){10}}
\put (180,84){\line(1,0){10}}
\put (196,100){\line(1,0){10}}

\put (145,103){\line(0,1){10}}
\put (145,87){\line(0,1){10}}
\put (145,71){\line(0,1){10}}
\put (145,55){\line(0,1){10}}
\put (145,39){\line(0,1){10}}
\put (161,103){\line(0,1){10}}
\put (161,87){\line(0,1){10}}
\put (161,71){\line(0,1){10}}
\put (161,55){\line(0,1){10}}
\put (177,103){\line(0,1){10}}
\put (177,87){\line(0,1){10}}
\put (177,71){\line(0,1){10}}
\put (193,103){\line(0,1){10}}
\put (193,87){\line(0,1){10}}
\put (209,103){\line(0,1){10}}

\put(124,98){1}
\put(124,82){2}
\put(124,66){3}
\put(124,50){4}
\put(124,34){5}

\put(142,117){1}
\put(158,117){2}
\put(174,117){3}
\put(190,117){4}
\put(206,117){5}

\put (142,84){\line(1,0){6}}
\put (145,81){\line(0,1){6}}

\put (142,100){\line(1,0){6}}
\put (145,97){\line(0,1){6}}

\put (142,52){\line(1,0){6}}
\put (145,49){\line(0,1){6}}

\put (158,100){\line(1,0){6}}
\put (161,97){\line(0,1){6}}

\put(205,104){\oval(8,8)[br]}
\put(189,88){\oval(8,8)[br]}
\put(173,72){\oval(8,8)[br]}
\put(157,56){\oval(8,8)[br]}
\put(141,40){\oval(8,8)[br]}

\put(189,104){\oval(8,8)[br]}
\put(197,96){\oval(8,8)[tl]}

\put(157,72){\oval(8,8)[br]}
\put(165,64){\oval(8,8)[tl]}

\put(141,72){\oval(8,8)[br]}
\put(149,64){\oval(8,8)[tl]}

\put(173,88){\oval(8,8)[br]}
\put(181,80){\oval(8,8)[tl]}

\put(157,88){\oval(8,8)[br]}
\put(165,80){\oval(8,8)[tl]}

\put(173,104){\oval(8,8)[br]}
\put(181,96){\oval(8,8)[tl]}
\end{picture}
\vspace{-36pt}

\caption{Two ways of drawing an rc-graph.}
\end{figure}

Denote by $\EuScript{RC}(w)$ the set of all rc-graphs which correspond to
$w$. For $D\in \EuScript{RC}(w)$ let $x^D=\prod_{(\alpha_k,
a_k-\alpha_k+1)\in D} x_{\alpha_k}.$ Now we state a theorem proved in
\cite {bjs} and in \cite {fs}.

\begin{theorem} For any permutation $w \in S_{\infty}$,
\begin{displaymath}
P_w = \sum_{{\bf a} \in R(w)}\, \sum_{\alpha_1...\alpha_p \in C({\bf a})}
x_{\alpha_1}...x_{\alpha_p} =\sum_{D\in \EuScript{RC}(w)}x^D.
\end{displaymath}
\end{theorem}
On an rc-graph we can define a ladder move (of size $m$) as shown in 
Figure 2. Note that after the move we get an rc-graph corresponding to the
same permutation (since the number of intersections is preserved
and the same strands intersect in the two figures).

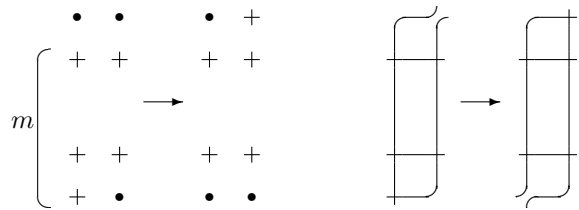
\begin{figure}
\begin{picture}(200, 130)

\put (22,56){\line(1,0){6}}
\put (25,53){\line(0,1){6}}
\put (38,56){\line(1,0){6}}
\put (41,53){\line(0,1){6}}
\put (22,40){\line(1,0){6}}
\put (25,37){\line(0,1){6}}
\put (22,92){\line(1,0){6}}
\put (25,89){\line(0,1){6}}
\put (38,92){\line(1,0){6}}
\put (41,89){\line(0,1){6}}

\put (25,108){\circle*{3}}
\put (41,108){\circle*{3}}
\put (41,40){\circle*{3}}

\put(0,66){$m$}
\put(15,66){\oval(10,60)[l]}
\put (50,76){\vector(1,0){15}}

\put (72,56){\line(1,0){6}}
\put (75,53){\line(0,1){6}}
\put (88,56){\line(1,0){6}}
\put (91,53){\line(0,1){6}}
\put (72,92){\line(1,0){6}}
\put (75,89){\line(0,1){6}}
\put (88,92){\line(1,0){6}}
\put (91,89){\line(0,1){6}}
\put (88,108){\line(1,0){6}}
\put (91,105){\line(0,1){6}}

\put (75,108){\circle*{3}}
\put (75,40){\circle*{3}}
\put (91,40){\circle*{3}}

\put (148,40){\line(1,0){10}}
\put (148,56){\line(1,0){10}}
\put (148,92){\line(1,0){10}}
\put (148,108){\line(1,0){10}}

\put (161,43){\line(0,1){10}}
\put (145,43){\line(0,1){10}}
\put (161,59){\line(0,1){30}}
\put (145,59){\line(0,1){30}}
\put (161,95){\line(0,1){10}}
\put (145,95){\line(0,1){10}}

\put (142,56){\line(1,0){6}}
\put (145,53){\line(0,1){6}}
\put (158,56){\line(1,0){6}}
\put (161,53){\line(0,1){6}}
\put (142,40){\line(1,0){6}}
\put (145,37){\line(0,1){6}}
\put (142,92){\line(1,0){6}}
\put (145,89){\line(0,1){6}}
\put (158,92){\line(1,0){6}}
\put (161,89){\line(0,1){6}}

\put(157,44){\oval(8,8)[br]}
\put(157,112){\oval(8,8)[br]}

\put(165,104){\oval(8,8)[tl]}
\put(149,104){\oval(8,8)[tl]}

\put (170,76){\vector(1,0){15}}

\put (198,40){\line(1,0){10}}
\put (198,56){\line(1,0){10}}
\put (198,92){\line(1,0){10}}
\put (198,108){\line(1,0){10}}

\put (211,43){\line(0,1){10}}
\put (195,43){\line(0,1){10}}
\put (211,59){\line(0,1){30}}
\put (195,59){\line(0,1){30}}
\put (211,95){\line(0,1){10}}
\put (195,95){\line(0,1){10}}

\put (192,56){\line(1,0){6}}
\put (195,53){\line(0,1){6}}
\put (208,56){\line(1,0){6}}
\put (211,53){\line(0,1){6}}
\put (192,92){\line(1,0){6}}
\put (195,89){\line(0,1){6}}
\put (208,92){\line(1,0){6}}
\put (211,89){\line(0,1){6}}
\put (208,108){\line(1,0){6}}
\put (211,105){\line(0,1){6}}

\put(207,44){\oval(8,8)[br]}
\put(191,44){\oval(8,8)[br]}

\put(199,36){\oval(8,8)[tl]}
\put(199,104){\oval(8,8)[tl]}
\end{picture}
\vspace{-36pt}
\caption{Ladder moves.}
\end{figure}

 We shall call the inverse operation an inverse ladder move. Billey and
Bergeron in \cite {bb} prove that all graphs in $\EuScript{RC}(w)$ can be
obtained by applying a succession of ladder moves to $D_{bot}(w) = \{(i,c)
: c \leq m_i\}$ (where $m_i = \#\{j : j>i \textrm{ and } w(j) > w(i)\}$),
which is the ``bottom'' rc-graph corresponding to $w$. Note that the
bottom graph has the property that the rows are left-justified. For
instance, the graph shown in Figure 1 is a bottom graph.

\section{Generalizing the insertion algorithm}

Recall that Pieri's rule seeks to give a formula for multiplication of
a Schubert polynomial by $P_{\sigma[r,m]}$ (a homogeneous symmetric
polynomial). In other words, it computes the generalized
Littlewood-Richardson coefficients $c^u_{vw}$ for $w=\sigma[r,m]$, where
the coefficients are defined by
\[
P_vP_w = \sum_u c^{u}_{vw} P_u.
\]
(We use the fact that Schubert polynomials form a basis of all
polynomials to get a unique such expansion.) To enable us to find this
formula, we need a modification of the insertion algorithm, which was
defined in \cite{bb} and used to compute a family of Littlewood-Richardson
coefficients in \cite{k}.

Notice that each rc-graph for the permutation $\sigma[r,m]$ can be given
by specifying the number of intersections in the rows on or above row
$r$. Assume this rc-graph is given by $k_s$ elements in row $i_s$,
where $r\geq i_1 > i_2 > ... > i_t> 0$ and $\sum k_s=m$. Below is the
algorithm which inserts this rc-graph into a general rc-graph.

\begin{algorithmn} 
Given an rc-graph $D$, and a level $r$, suppose we have
to insert $k_s$ elements into row $i_s$ of $D$, where $r\geq i_1 > i_2 >
... > i_t$. We keep an ordered sequence of $(a_i,b_i)$'s
during the course of the algorithm, such that $a_1\leq a_2\leq...\leq
a_l\leq r$, each $b_i>r$ and all $b_i$'s are distinct. Initially the sequence is empty.

Proceed row by row as follows. Starting with row $i_1$, find the rightmost
position $(i_1, j)$ where the configuration is one of the following:
\begin{figure}[!h]
\vspace{-30pt}
\begin{picture}(400, 100)

\put (50,55){\line(1,0){10}}
\put (63,42){\line(0,1){10}}
\put (66,55){\line(1,0){10}}
\put (63,58){\line(0,1){10}}

\put(59,59){\oval(8,8)[br]}
\put(67,51){\oval(8,8)[tl]}

\put(44,53){$s$}
\put(60,30){$q$}

\put(90,60){with $s\leq r <q$,}
\put(90,45){$q\neq b_j$ for any $j$.}

\put (200,55){\line(1,0){10}}
\put (213,42){\line(0,1){10}}
\put (216,55){\line(1,0){10}}
\put (213,58){\line(0,1){10}}

\put(209,59){\oval(8,8)[br]}
\put(217,51){\oval(8,8)[tl]}

\put(190,53){$b_i$}
\put(210,30){$q$}

\put(240,60){with $r <q$,}
\put(240,45){$q\neq b_j$ for any $j$.}
\end{picture}
\vspace{-36pt}
\caption{$ $}
\end{figure}

Add this intersection. If we are in the first case, insert $(s,q)$ into
the sequence $(a_i,b_i)$ in the rightmost position, such that $a_i$'s
remain nondecreasing in the sequence. 
($(s,q)$ are the rows where the two strands shown in  Figure 3 originate.)
If we are in the second case, add
$(a_i,q)$ just before where $(a_i,b_i)$ is in the sequence. Repeat the
above process $k_1$ times to insert $k_1$ intersections into the row~$i_1$
to get a new graph together with a sequence $(a_1,b_1),....,(a_{k_1},
b_{k_1})$ where we know for sure that the $b_i$'s are distinct and greater
than $r$ and the $a_i$'s are ordered and less than or equal to $r$. Also, the
graph we have obtained might not be an rc-graph. So we go up to the next
row and perform the following rectification procedure. Starting from the
left, we look for either of the two configurations in Figure 4. 
\begin{figure}[!h]
\vspace{-30pt}
\begin{picture}(400, 100)

\put (50,55){\line(1,0){26}}
\put (63,42){\line(0,1){26}}

\put(40,53){$b_i$}
\put(60,30){$a_i$}

\put (200,55){\line(1,0){26}}
\put (213,42){\line(0,1){26}}

\put(190,53){$b_i$}
\put(210,30){$b_j$}

\put(240,53){with $a_i=a_j$,}
\end{picture}
\vspace{-36pt}
\caption{$ $}
\end{figure}

We remove the intersection, and delete the pair $(a_i,b_i)$ from the sequence
of $(a,b)$'s. If the intersection was at position $(i,j)$, find the
maximal $j' < j$ such that the configuration at $(i,j')$ is of the form
shown in Figure 3. (Such a position must exist, since the $b_i$ is greater
than $r$ and $(a_i,b_i)$ is no longer in the sequence.) We add the
intersection, and to the list of $(a,b)$'s we add the new pair the same
way it was described before. Then we look to the right to see if there is
another intersection of the form given in Figure 4. If it exists, remove
it together with the appropriate $(a_i,b_i)$ and proceed as before. We call
this process rectification of a given row.

When we are done with rectifying the current row, we add what
intersections we must to this row in the same manner as for the previous
row, giving us some new pairs $(a_j,b_j)$. We then proceed to the next row
and repeat the procedure. It is clear from the construction that the $b$'s
remain distinct and the $a$'s remain ordered. 
\qed
\end{algorithmn}

The rest of this section is concerned with proving the following:

\begin{theorem}
\label{alg}
The above algorithm produces an rc-graph.
\end{theorem}

Before proving Theorem \ref{alg}, we will need two lemmas and another
algorithm.

\begin{lemma}
\label{order} After row $\ell$ has been rectified, the strands $b_i,
..., b_{i+k}$ with  $a=a_i=...=a_{i+k}$ pass from row $\ell$ to 
row $\ell -1$  to the  left of the place where the strand $a$ passes from
row $\ell$ to row $\ell-1$. Moreover, they pass from row
$\ell$ to row $\ell-1$ in  the same  order they appear in the $(a,b)$
sequence. (See Figure 5 for an illustration.)
\end{lemma}

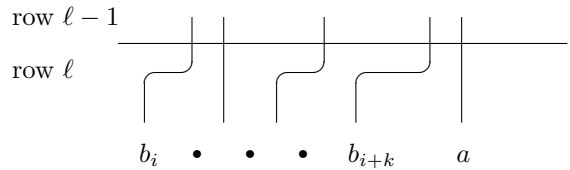
\begin{figure}[!h]
\vspace{-12pt}
\begin{picture}(220, 100)

\put (10,62){row $\ell$}
\put (10,82){row $\ell-1$}
\put (50,75){\line(1,0){170}}

\put (60,45){\line(0,1){15}}
\put (90,45){\line(0,1){40}}
\put (110,45){\line(0,1){15}}
\put (140,45){\line(0,1){15}}
\put (180,45){\line(0,1){40}}

\put(64,60){\oval(8,8)[tl]}
\put(114,60){\oval(8,8)[tl]}
\put(144,60){\oval(8,8)[tl]}

\put(74,68){\oval(8,8)[br]}
\put(124,68){\oval(8,8)[br]}
\put(164,68){\oval(8,8)[br]}

\put (64,64){\line(1,0){10}}
\put (114,64){\line(1,0){10}}
\put (144,64){\line(1,0){20}}

\put(78,68){\line(0,1){17}}
\put(128,68){\line(0,1){17}}
\put(168,68){\line(0,1){17}}

\put (58,30){$b_i$}
\put (137,30){$b_{i+k}$}
\put (178,30){$a$}

\put (80,33) {\circle*{3}}
\put (100,33) {\circle*{3}}
\put (120,33) {\circle*{3}}
\end{picture}
\vspace{-24pt}
\caption{The order in which the strands
$b_i,...,b_{i+k}$ and $a$ intersect the line, which
separates rows $\ell$
and $\ell -1$, is the same as the order of $b_i$'s in the $(a,b)$
sequence.}
\end{figure}

\begin{proof} This lemma can be proved by induction on $\ell$. At the
beginning of the algorithm, set $\ell= r-1$. Then the statement of the
lemma holds, since the $(a,b)$ sequence is empty.

Assume we know that the lemma holds for $\ell+1$. Notice that if during the
rectification of row $\ell$ no pairs $(a_j,b_j)$ with $a_j=a$ were
added or removed from the $(a,b)$ sequence, the strands $b_i,...,b_{i+k}$
and $a$ never intersect each other in row $\ell$. Hence they
pass between rows $\ell$ and $\ell-1$ in the same order they pass from
row $\ell+1$ to row $\ell$.

Moreover, the order of the strands $b_i,...,b_{i+k}$ and $a$ can change
only when some of them intersect in row $\ell$. But then this
intersection is one of the intersections in Figure 4. The algorithm removes
this intersection, and it removes the appropriate tuple from the $(a,b)$
sequence to preserve the order of the stands.

At the same time, when an intersection is added, the addition of
$(a_j,b_j)$ to the list is consistent with order of the strands
$b_i,...,b_{i+k}$ and $a$ passing from row $\ell$ to row $\ell-1$.
\end{proof}

\begin{lemma}
\label{ab}
The permutation carried by the strands of the graph at each moment of the
algorithm is $wt_{a_1b_1}...t_{a_{m'}b_{m'}}$ where
$\langle(a_1,b_1),...,(a_{m'},b_{m'})\rangle $ 
is the sequence up to that point in the
algorithm.
\end{lemma}

\begin{proof} We again use induction: we need to show that the above statement
remains true when we insert an intersection or delete one.

Assume we are in the first case of Figure $3$; that is, we are adding an
intersection and $(s,q)$ to the list of $(a,b)$'s. We can easily see that
adding this intersection multiplies the permutation of the graph by
$t_{s,q}$ on the right. So, if $\langle (a_1,b_1),...,(a_{m'},b_{m'})\rangle$ is the
list before this insertion, the permutation of the graph after the
insertion is $wt_{a_1b_1}...t_{a_{m'}b_{m'}}t_{s,q}$. Using the fact that
$t_{ab}$ and $t_{cd}$ commute when $a,b,c,d$ are all distinct, we can
commute $t_{s,q}$ trough $t_{a_1b_1}...t_{a_{m'}b_{m'}}$ to the left until
some $a_j\leq s$. Then the permutation carried by our graph is
$wt_{a_1b_1}...t_{a_jb_j}t_{sq}...t_{a_{m'}b_{m'}}$. So, after inserting
$(s,q)$ into the list after $(a_i,b_j)$, the property we are trying to
prove holds.

If we are in the second situation of Figure 3, then we multiply the
permutation $wt_{a_1b_1}...t_{a_{m'}b_{m'}}$ by $t_{b_i,q}$ on the
right. Then $t_{b_i,q}$ can be moved to the left until the place $i$,
since all $b_j$'s are distinct and thus $t_{b_i,q}$ commutes with all
$t_{a_j,b_j}$ if $j>i$. After that we use the following identity, which
shows why we have to add $(a_i,q)$ to the $(a,b)$ list before $(a_i,b_i)$:
\begin{equation}
\label{commutation}
t_{a,b_i}t_{b_i,q} = t_{a,q}t_{a,b_i}.
\end{equation}

For the deletion, we use almost identical arguments in both cases
of Figure 4, since the deletion of an intersection also multiplies
the permutation of the graph by the appropriate transposition. The only
difference is that in the second case, instead of (\ref{commutation}) we
have to use the following identity:
\[
t_{a,b_i}t_{a,b_j}t_{b_i,b_j} = t_{a,b_j}.
\]
This concludes the proof of the lemma.
\end{proof}

To give a proof of Theorem \ref{alg} we have to introduce another
algorithm (call it Algorithm 2) that we perform on the intermediate graph.
This algorithm will take the graph which was rectified upto the row
$\ell$, and remove some intersections to produce an rc-graph of the
original permutation $w$.

\addtocounter{algorithm}{1}
\begin{algorithm} We start with the list of $(a_i,b_i)$'s and, starting
with the last row $\ell$ we have finished with, go down row by row and
from right to left in each row. We look for intersections of the form
\begin{figure}[!h]
\vspace{-24pt}
\begin{picture}(400, 100)

\put (50,55){\line(1,0){26}}
\put (63,42){\line(0,1){26}}

\put(40,53){$a_i$}
\put(60,30){$b_i$}

\put (200,55){\line(1,0){26}}
\put (213,42){\line(0,1){26}}

\put(190,53){$b_j$}
\put(210,30){$b_i$}

\put(240,53){with $a_i=a_j$,}

\end{picture}
\vspace{-36pt}
\caption{$ $}
\end{figure}

Whenever we see these, we remove the intersection and $(a_i,b_i)$ from its
place in the list. Note that Lemma \ref{order} and
Lemma \ref{ab} hold at each moment of this algorithm. Ultimately all the
$(a_i,b_i)$ will be removed from the list because each strand $a_i$
definitely intersects the strand $b_i$ at each moment of Algorithm 2 by
Lemma \ref{order}, so all of them will get removed and at the end we will
get a graph and an empty list of $(a,b)$'s.
\qed
\end{algorithm} 

The result of Algorithm 2 will be a graph with the permutation $w$ by
Lemma \ref{ab}, with exactly $l(w)$ intersections. Therefore it will be an
rc-graph for $w$. We will use this fact to prove that the main algorithm
generates an rc-graph.

We come back to the proof of Theorem \ref{alg}. Suppose we have rectified
some row $\ell$ and inserted some elements in it, and we want to show
there are no double intersections below (and including) that row. Because of
the inductive hypothesis, this works for the previous row $\ell+1$ (the base
case is clear). Hence we need to check that no two strands intersect at
row $\ell$ and at some lower row after the algorithm has rectified 
row $\ell$ and inserted all the required intersections into it.

We will show that the only double intersection which might be introduced
in row $\ell$ by insertions in the lower rows will be removed by the
algorithm and no new double intersections will be introduced in this row.

Let us show that after rectifying row $\ell+1$, and deleting all the
intersections from Figure 4 in row $\ell$ at once, there are no
strands which intersect twice below row $\ell$. For this apply
Algorithm 2 to the graph constructed after rectifying row $\ell+1$.
Then the resulting graph is an rc-graph $D'$ of the permutation $w$. Start
applying the inverse Algorithm 2 to this rc-graph $D'$; that is, start
adding intersections to $D'$ in the order opposite to the order they were
removed.  It is then not difficult to see that since we are only adding
intersections from Figure 6 and at each moment there are no double
intersections below row $\ell+1$, the only double intersections at row
$\ell$ have to look like the ones shown in Figure 4. This shows that
by removing intersections from row $\ell$ we remove all the double
intersections below or at this row.

Now let's check that we do not create any double intersections below or at
row $\ell$ during rectification of this row or during the insertion of new
intersections into it. Clearly the only way we could have done this is if
we inserted an intersection in the position as in the second case of
Figure 3, so that the strands $b_i$ and $q$ intersected before, in a lower
row. (In the first case of Figure 3, the two strands could not have
intersected in a lower row, since $s <q$.) So, we shall show that in the
second case the strands $b_i$ and $q$ cannot intersect in a lower row, that
is, $b_i < q$. Proving this will finish the proof of Theorem \ref{alg}.

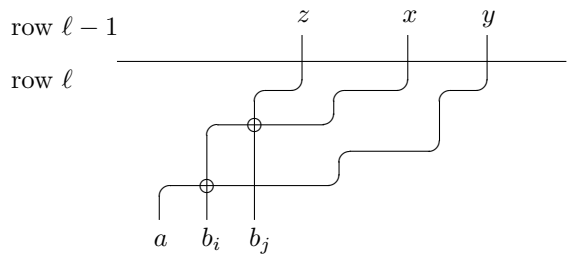
\begin{figure}[b] 
\vspace{-18pt}
\begin{picture}(220, 130)

\put (10,115){row $\ell-1$}
\put (10,95){row $\ell$}
\put (50,105){\line(1,0){170}}

\put(106,90){\oval(8,8)[tl]}
\put(136,90){\oval(8,8)[tl]}
\put(176,90){\oval(8,8)[tl]}

\put(116,98){\oval(8,8)[br]}
\put(156,98){\oval(8,8)[br]}
\put(186,98){\oval(8,8)[br]}

\put(120,98){\line(0,1){17}}
\put(160,98){\line(0,1){17}}
\put(190,98){\line(0,1){17}}

\put(106,94){\line(1,0){10}}
\put(136,94){\line(1,0){20}}
\put(176,94){\line(1,0){10}}

\put(102,45){\line(0,1){45}}
\put(132,85){\line(0,1){5}}
\put(172,75){\line(0,1){15}}

\put(128,85){\oval(8,8)[br]}
\put(168,75){\oval(8,8)[br]}

\put(88,81){\line(1,0){40}}
\put(138,71){\line(1,0){30}}

\put(88,77){\oval(8,8)[tl]}
\put(138,67){\oval(8,8)[tl]}

\put(84,45){\line(0,1){32}}
\put(134,62){\line(0,1){5}}

\put(130,62){\oval(8,8)[br]}

\put(70,58){\line(1,0){60}}

\put(70,54){\oval(8,8)[tl]}

\put(66,45){\line(0,1){9}}

\put (84,58){\circle{5}}
\put (102,81){\circle{5}}

\put(64,35){$a$}
\put(82,35){$b_i$}
\put(100,35){$b_j$}

\put(118,120){$z$}
\put(158,120){$x$}
\put(188,120){$y$}
\end{picture}
\vspace{-30pt}
\caption{During Algorithm 2, $b_i$ moved from $x$ to $z$
and then to $y$ in row $\ell$.}
\end{figure}

After row $\ell$ is rectified, we again apply Algorithm 2. Notice that
this algorithm keeps moving strand $b_i$ to the left in row $\ell$
till the last point, where $(a_i,b_i)$ is removed from the list of
$(a,b)$'s, at which point strand $b_i$ moves right in row $\ell$,
even to the right of its {\em original} position. For example, in Figure
7, the circled intersections are removed during Algorithm 2. When the
first one is removed, $(a_j,b_j)$ is removed from the list and $b_i$ moves
from $x$ to the left, to $z$, and when the second intersection is removed,
$(a_i,b_i)$ is removed and $b_i$ moves to the right of $x$, to $y$.

Therefore the position $x$ of strand $b_i$ at row $\ell$ at the
start of Algorithm 2 is to the left of the position $y$ at row $\ell$
at the end of it.

Now we look at what happens after $(a_i,b_i)$ is moved off the list in
Algorithm 2, by removing an intersection with $b_{i_1}$ (it moved first 
from $x$ to $z$ and then from $z$ to $y$ in row $\ell$). We look next
at the history of $b_{i_1}$ till $(a_{i_1},b_{i_1})$ is moved off the list,
and so on, getting a figure like Figure 8. 

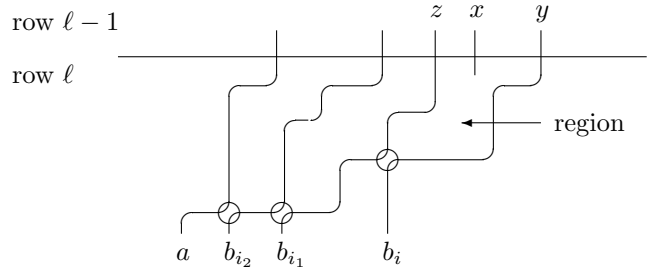
\begin{figure} 
\begin{picture}(220, 130)

\put(110,98){\line(0,1){17}}
\put(150,98){\line(0,1){17}}
\put(170,88){\line(0,1){27}}
\put(185,98){\line(0,1){17}}
\put(210,98){\line(0,1){17}}

\put(106,98){\oval(8,8)[br]}
\put(146,98){\oval(8,8)[br]}
\put(166,88){\oval(8,8)[br]}
\put(206,98){\oval(8,8)[br]}

\put(96,94){\line(1,0){10}}
\put(131,94){\line(1,0){15}}
\put(156,84){\line(1,0){10}}
\put(196,94){\line(1,0){10}}

\put(96,90){\oval(8,8)[tl]}
\put(131,90){\oval(8,8)[tl]}
\put(156,80){\oval(8,8)[tl]}
\put(196,90){\oval(8,8)[tl]}

\put(92,50){\line(0,1){40}}
\put(127,85){\line(0,1){5}}
\put(152,70){\line(0,1){10}}
\put(192,70){\line(0,1){20}}

\put(88,50){\oval(8,8)[br]}
\put(123,85){\oval(8,8)[br]}
\put(148,70){\oval(8,8)[br]}
\put(188,70){\oval(8,8)[br]}

\put(78,46){\line(1,0){10}}
\put(117,81){\line(1,0){5}}
\put(138,66){\line(1,0){10}}
\put(156,66){\line(1,0){32}}

\put(78,42){\oval(8,8)[tl]}
\put(117,77){\oval(8,8)[tl]}
\put(138,62){\oval(8,8)[tl]}
\put(156,62){\oval(8,8)[tl]}

\put(74,38){\line(0,1){4}}
\put(113,50){\line(0,1){27}}
\put(134,50){\line(0,1){12}}
\put(152,38){\line(0,1){24}}

\put(109,50){\oval(8,8)[br]}
\put(130,50){\oval(8,8)[br]}

\put(96,46){\line(1,0){13}}
\put(116,46){\line(1,0){14}}

\put(96,42){\oval(8,8)[tl]}
\put(116,42){\oval(8,8)[tl]}

\put(92,38){\line(0,1){4}}
\put(112,38){\line(0,1){4}}

\put (92,46){\circle{8}}
\put (112,46){\circle{8}}
\put (152,66){\circle{8}}

\put(72, 28){$a$}
\put(90,28){$b_{i_2}$}
\put(110,28){$b_{i_1}$}
\put(150,28){$b_i$}

\put(168,120){$z$}
\put(183,120){$x$}
\put(208,120){$y$}

\put (210,80){\vector(-1,0){30}}
\put (215, 77){region}

\put (10,115){row $\ell-1$}
\put (10,95){row $\ell$}
\put (50,105){\line(1,0){200}}

\end{picture}
\vspace{-24pt}
\caption{No strands between $a$ and $b_i$ can pass
in the shown region.}
\end{figure}

Note that no two consecutive strands shown in the figure can intersect
above the encircled spot (where one of them is moved off the list),
because that intersection would then have to be removed earlier than the
encircled one, by the construction of Algorithm 2. Therefore, from the
fact that the outcome of Algorithm 2 is an rc-graph, we conclude that in
the region shown, there can pass no strand between $a$ and $b_i$
(otherwise some double intersections must occur). Hence the $q$ from the
second case of Figure 3 must always be greater than $b_i$. This concludes
the proof of Theorem \ref{alg}.
\qed

\begin{remark}
Assume that the permutation $w$ satisfies the 
condition $w(i+1)>w(i)$ if $i>r$. Each rc-graph of such a permutation has
the following defining property: no two strands which start below row $r$
intersect. In this case it can be shown that the above algorithm is
equivalent to one by one insertion. This is consistent with the results in
\cite{k}, where the Pieri formula in this special case was proved using
one by one insertion. The above statement and the fact that
one by one insertion does not work in general indicate that this
algorithm provides the proper generalization of the Bergeron-Billey algorithm
of \cite{bb}.
\end{remark}

\section{The proof of Pieri's formula}
To prove Pieri's formula using the above insertion algorithm, we will
construct an inverse algorithm, which takes an rc-graph for $w' =
wt_{a_1,b_1}...t_{a_m,b_m}$, such that the $b_i$'s are distinct and
greater than $r$ and $a_1 \leq a_2 \leq ... \leq a_m\leq r$, satisfying
$l(w')= l(w)+m$ and $w'(b_i) < w'(b_j)<w'(a_i)$ for every $i<j$ such that
$a_i=a_j$, and produces an rc-graph for $w$ and an rc-graph for
$\sigma[r,n]$ (remember that rc-graphs of $\sigma[r,m]$ are given by
specifying the number of intersections $k_s$ in each row $s\leq r$ 
such that $\sum k_s=m$). Given the insertion algorithm and its inverse, we
will clearly produce a bijection:
\[
\EuScript{RC}(w) \times \{(k_1,...,k_r) : \sum k_i = m\}
=\EuScript{RC}(w) \times \EuScript{RC}(\sigma[r,n])
\rightarrow \bigcup \EuScript{Rc}(w'),
\]
where $w'$ ranges as stated above.

Below we will produce this inverse algorithm, which will automatically
prove the following theorem. 

\begin{theorem} 
\label{pieri}
\[
P_w P_{\sigma[r,m]} = \sum P_{w'},
\]
where the sum is over all $w' = wt_{a_1,b_1}...t_{a_m,b_m}$ such that the
$b_i$'s are distinct and greater than $r$, $a_1 \leq a_2 \leq ... \leq
a_m\leq r$, satisfying $l(w')= l(w)+m$ and $w'(b_i) < w'(b_j)<
w'(a_i)$ for every $i<j$ with $a_i=a_j$.
\end{theorem}

\begin{inversealgorithm}

We are given an rc-graph $R'$ for a permutation 
\[w' =
wt_{a_1b_1}...t_{a_mb_m}\]
with distinct $b_i$'s which are greater
than $r$ and with $a_1 \leq a_2 ...\leq a_m\leq r$, $l(w') = l(w) +
m$, and $w'(a) > w'(b) > w'(c)$ for all $a,b,c$ such that $(a,b)$ is
some $(a_s,b_s)$ and $(a,c)$ is some $(a_t,b_t)$ with $s > t$. The
algorithm is defined as follows: starting from the top (first) row of
the rc-graph, we look from right to left for the occurrence of one of
the configurations of Figure 9. 

\begin{figure}
\vspace{-24pt}
\begin{picture}(400, 100)

\put (50,55){\line(1,0){26}}
\put (63,42){\line(0,1){26}}

\put(40,53){$a_i$}
\put(60,30){$b_i$}

\put (200,55){\line(1,0){26}}
\put (213,42){\line(0,1){26}}

\put(190,53){$b_j$}
\put(210,30){$b_i$}

\put(240,53){with $a_i=a_j$,}

\end{picture}
\vspace{-36pt}
\caption{$ $}
\end{figure}

We then remove this intersection from the rc-graph and consequently
remove $(a_i,b_i)$ from the list of $(a,b)$'s.  Immediately after we
remove such an intersection we look to its right for a configuration
of the form shown in Figure 10. 

\begin{figure}
\vspace{-24pt}
\begin{picture}(400, 100)

\put (50,55){\line(1,0){10}}
\put (63,42){\line(0,1){10}}
\put (66,55){\line(1,0){10}}
\put (63,58){\line(0,1){10}}

\put(59,59){\oval(8,8)[br]}
\put(67,51){\oval(8,8)[tl]}

\put(44,53){$q$}
\put(60,30){$s$}

\put(90,60){with $s\leq r <q$,}
\put(90,45){$q\neq b_i$ for any $i$.}

\put (200,55){\line(1,0){10}}
\put (213,42){\line(0,1){10}}
\put (216,55){\line(1,0){10}}
\put (213,58){\line(0,1){10}}

\put(209,59){\oval(8,8)[br]}
\put(217,51){\oval(8,8)[tl]}

\put(190,53){$q$}
\put(210,30){$b_i$}

\put(240,60){with $r <q$,}
\put(240,45){$q\neq b_j$ for any $j$.}

\end{picture}
\vspace{-36pt}
\caption{$ $}
\end{figure}

If we do not find such a configuration, we say that an intersection is
removed from this row and we move on with the algorithm, looking left
along that row (and after that, going to the next row) for an intersection
of the form shown in Figure 9. (At the same time we record the number of removed
intersections from this row, as these numbers will give an rc-graph for
$\sigma[r,m]$ at the end of the algorithm.) If we find such a
configuration, we add the intersection at this place (i.e. the first
configuration of this sort to the right of the removed
intersection) and add $(q,s)$ or $(a_i,q)$ to the list of $(a,b)$'s the
same way it was done in the forward insertion algorithm, i.e we add
$(q,s)$ to the $(a_i, b_i)$'s to preserve the order of $a$'s ($a_{i-1}\leq
a_i< a_{i+1}$) and we add $(a_i,q)$ right after $(a_j,b_j)$. Then we
proceed with the algorithm, looking left for more intersections to remove.
We do this till the list of $(a,b)$'s becomes empty. At the end we are
left with a graph for the permutation $w$. Moreover, the previously
recorded numbers of removed intersections in each row where intersections
were removed from the rc-graph produce an rc-graph for $\sigma[r,m]$. \qed
\end{inversealgorithm}

The fact that the resulting graph is an rc-graph can be shown using
the same sort of technique used for the forward algorithm. It is also
evident that the algorithm is the inverse of the forward algorithm.
This finishes the proof of Theorem \ref{pieri}.

Let us now restate Pieri's formula in its original form, as it appeared in
\cite{ls1}. Assume that for the $(a,b)$ sequence from Theorem \ref{pieri}
we have $a_1=...=a_{i_1}< a_{i_1+1}=... =a_{i_2}<...<a_{i_p}=...=a_m$.
Then $w'=w\zeta_1...\zeta_p$, where each
$\zeta_k=(a_{i_{k-1}+1},b_{i_{k-1}}+1)...(a_{i_k},b_{i_k})$ is a cycle.
We will say that the cycle $\zeta_k$ has size $s(\zeta_k)=i_k-i_{k-1}$.
Additionally, it follows from Theorem \ref{pieri} that for each $\zeta_k$,
there exists exactly one $a$ with $w(a)<w\zeta_k(a)$ and exactly one
$a'\leq r$ with $w(a')\neq w\zeta_k(a')$. (Actually, $a=a'=a_{i_k}$.)
Hence it is easy to conclude that Theorem \ref{pieri} is equivalent to the
following formulation of Pieri's formula given in \cite{ls1}.

\begin{theorem}
\[
P_w P_{\sigma[r,m]} = \sum P_{w'},
\]
where the sum is over all $w' = w\zeta_1...\zeta_p$ satisfying the
following conditions: $l(w')-l(w)=m=\sum s(\zeta_k)$, and each $\zeta_k$ is a
cycle such that  there exists exactly one $a$ with $w(a)<w\zeta_k(a)$ and
exactly one $a'\leq r$ with $w(a')\neq w\zeta_k(a')$.
\end{theorem}

\end{document}